\documentclass[11pt]{article}
    
\usepackage{amsmath,amsfonts,amsthm,amscd,amssymb,graphicx}

\numberwithin{equation}{section}

\usepackage{hyperref}
\usepackage{color}

\newtheorem{theorem}{Theorem}[section]
\newtheorem{theo}{Theorem}[section]

\newtheorem{lem}[theorem]{Lemma}

\newtheorem{prop}[theorem]{Proposition}

\def\eps{\varepsilon }

\newcommand{\RR}{\mathbb{R}}

\newcommand{\CC}{\mathbb{C}}

\newcommand{\NN}{{\mathbb N}}

\newcommand{\TT}{{\mathbb T}}

\def\beq{\begin{equation}}
\def\eeq{\end{equation}}
\def\bb1{{1\!\!1}}
\def\rit{{\Bbb R}}
\def\cit{{\Bbb C}}

\def\tit{{\Bbb T}}

\def\eps{\varepsilon}

\def\OS{\mathrm{OS}}

\begin{document}

\title{Green function of Orr Sommerfeld equations away from critical layers} 

\author{Emmanuel Grenier\footnotemark[1]
  \and Toan T. Nguyen\footnotemark[2]
}

\maketitle

\renewcommand{\thefootnote}{\fnsymbol{footnote}}

\footnotetext[1]{Equipe Projet Inria NUMED,
 INRIA Rh\^one Alpes, Unit\'e de Math\'ematiques Pures et Appliqu\'ees., 
 UMR 5669, CNRS et \'Ecole Normale Sup\'erieure de Lyon,
               46, all\'ee d'Italie, 69364 Lyon Cedex 07, France. Email: Emmanuel.Grenier@ens-lyon.fr}

\footnotetext[2]{Department of Mathematics, Penn State University, State College, PA 16803. Email: nguyen@math.psu.edu. 
TN's research was supported by the NSF under grant DMS-1405728. Part of this work was completed when TN was visiting the ENS de Lyon.}

\begin{abstract}
The classical Orr-Sommerfeld equations are the resolvent equations of 
the linearized Navier Stokes equations around a stationary shear layer profile in the half plane. 
In this paper, we derive pointwise bounds on the Green function of the Orr Sommerfeld problem 
away from its critical layers. 
\end{abstract}


\section{Introduction}


In this paper, we are interested in the study of linearized Navier Stokes equations around a given fixed profile $U_s = (U(z),0)$
as the viscosity goes to $0$. Namely, we consider the following set of equations
\beq \label{Nlin1}
\partial_t v + U_s \cdot \nabla v + v\cdot  \nabla U_s + \nabla p - \nu \Delta v = F,
\eeq
\beq \label{Nlin2}
\nabla\cdot  v = 0,
\eeq
on the half plane $x \in \TT$, $z \ge 0$,
with Dirichlet boundary condition
\beq \label{Nlin3}
v = 0 \quad \hbox{on} \quad z = 0 .
\eeq
We focus on the periodic case $x \in \tit$, the whole line case $x \in \rit$ being similar.
Throughout this paper, the background profile $U(z)$ is assumed to be sufficiently smooth, to satisfy $U(0)=0$ and   
\begin{equation}\label{def-Ubl} 
 |\partial_z^k (U(z) - U_+)| \le C_k e^{-\eta_0 z} , \qquad \forall~ z\ge 0, \quad k\ge 0,
 \end{equation}
for some finite constant $U_+$ and some positive constants $C_k$ and $\eta_0$. 

The inviscid limit problem \eqref{Nlin1}-\eqref{def-Ubl}  is a very classical problem 
that has led to a huge physical and mathematical literature, focussing in particular on the linear stability, 
on the dispersion relation, on the study of eigenvalues and eigenmodes, and on the onset of nonlinear instabilities and turbulence
(see \cite{Reid} for an introduction on these
topics, and the classical achievements of Rayleigh, Orr, Sommerfeld, Heisenberg, Tollmien, C.C. Lin,  and Schlichting).

Two cases arise. Either the profile $U$ is linearly stable for the corresponding linearized Euler equations (the case when $\nu = 0$)
 or it is linearly unstable for these limiting equations.
 In this paper, we consider the unstable case,
  leaving the stable case to be treated in \cite{GrN4}, which turns out to be much delicate.  
In the unstable case, it is well known \cite{Grenier00CPAM} that the profile $U$ is linearly unstable for the
linearized Navier Stokes equations provided $\nu$ is sufficiently
small, or equivalently, the Reynolds number $R = \nu^{-1}$ is sufficiently large. 
However, in order to go from linear to nonlinear instability, more precise information on solutions to the linearized problem is required. 
Let us mention several efforts in treating the stability and instability of nonlinear boundary layers in the small viscosity limit 
\cite{WE,Rousset,DGV, GVMM,Grenier00CPAM,GuoN,Mae, SC2}.

A natural and traditional approach to study linearized Navier Stokes equations is to take the Fourier Laplace transform of these equations.
For this, in order to take advantage of the incompressibility relation (\ref{Nlin2}),
we introduce the stream function $\psi$ of $v$, defined by
$$
v = \nabla^\perp \psi ,
$$
and take its Fourier transform in the $x$ variable, with wave number $\alpha$, and Laplace transform in time, with Laplace variable
$$
\lambda = - i \alpha c,
$$
following historical notations. 
Equivalently, we focus on solutions $v$ of linearized Navier Stokes equations of the form
$$
 v  = \nabla^\perp \Big( e^{i\alpha (x-c t)} \phi(z)\Big),
$$ 
with source term of the same form. This leads to the classical Orr-Sommerfeld equation
\begin{equation}\label{OS1a}
(U-c) (\partial_z^2 - \alpha^2) \phi - U'' \phi = \epsilon (\partial_z^2 - \alpha^2)^2 \phi 
-i \alpha^{-1} f , \qquad \epsilon = \frac{\nu}{i\alpha } 
\end{equation} 
on the half line $z\ge 0$, 
together with the boundary conditions
\begin{equation}\label{OS2}
\phi_{\vert_{z=0}} = \phi'_{\vert_{z=0}} = 0, \qquad \lim_{z\to \infty} \phi(z) =0.
\end{equation}
Here, $\alpha \in \NN^*= \NN\setminus\{0\}$ denotes the tangential wave number and  $c \in \CC$ is the complex phase velocity. 

For the mathematical analysis, it is more convenient to multiply (\ref{OS1a}) by $i \alpha$, which leads to
\beq \label{OS1}
(\lambda + i\alpha U ) (\partial_z^2 - \alpha^2) \phi - i\alpha U'' \phi - \nu (\partial_z^2 - \alpha^2)^2 \phi = f .
\eeq
Such a spectral formulation of the linearized Navier-Stokes equations near a boundary layer shear profile has been 
intensively studied  in the physical literature.  We in particular refer to
\cite{Reid, LinBook, Sch} for the major works of Heisenberg, Tollmien, C.C. Lin, and Schlichting on the subject. 
We also refer to \cite{GGN1, GGN3, GGN2} for the rigorous spectral analysis of the Orr-Sommerfeld equations.  

In this paper,  we shall derive pointwise bounds on the Green function of the Orr-Sommerfeld problem \eqref{OS2}-\eqref{OS1}. For convenience, let us denote 
$$
 \Delta_\alpha: = \partial_z^2 - \alpha^2 
 $$
and
\begin{equation}\label{def-OS}
\OS_{\alpha,\lambda}(\phi): = (\lambda + i\alpha U) \Delta_\alpha \phi - i\alpha U'' \phi  -  \nu \Delta_\alpha^2 \phi .
\end{equation}
For each fixed $\alpha \in \NN^*$ and $\lambda\in \CC$, we denote by $G_{\alpha,\lambda}(x,z)$ the corresponding Green kernel 
of the Orr-Sommerfeld problem. 
By definition, for each $x\in \RR_+$, $G_{\alpha,\lambda}(x,z)$ solves 
$$ 
\OS_{\alpha,\lambda}(G_{\alpha,\lambda} (x,\cdot)) = \delta_x (\cdot)
$$
on $z\ge 0$, together with the boundary conditions:
\begin{equation}\label{G-noslip}
G_{\alpha,\lambda}(x,0) = \partial_z G_{\alpha,\lambda} (x,0) =0, \qquad \lim_{z\to \infty} G_{\alpha,\lambda}(x,z) =0.
\end{equation}
The Green kernel allows to solve the inhomogenous Orr-Sommerfeld problem 
\begin{equation}\label{OS-intro} 
\OS_{\alpha,\lambda} (\phi) = f,
\end{equation}
 or equivalently the resolvent equations of the linearized Navier-Stokes operator,  
through the following explicit expression for the solution $\phi$ 
$$ 
\phi(z) = \int_0^\infty G_{\alpha,\lambda} (x,z) f(x) \; dx .
$$
To construct the Green function, let us first note that as $z \to +\infty$ the homogenous Orr-Sommerfeld equation "converges" 
to  the following constant-coefficient equation
\begin{equation}\label{OS-plus} 
\OS_+(\phi) = (\lambda + i \alpha U_+ ) \Delta_\alpha \phi  -  \nu \Delta_\alpha^2 \phi =0,
\end{equation}
where $U_+ = \lim_{z\to \infty}U(z)$. This constant-coefficient equation has four independent solutions $e^{\pm\mu_s z}$ and $e^{\pm \mu_f^+ z}$, with 
\begin{equation}\label{def-sfrate}
\mu_s = | \alpha |, \quad
 \mu_f (z)=  \nu^{-1/2} \sqrt{ \lambda + \nu \alpha^2+i \alpha U(z) },\quad
  \mu_f^+  = \lim_{z\to \infty}\mu_f(z),
 \end{equation}
in which we take the positive real part of the square root. 

As will be proved later, there exist four solutions to the homogenous Orr-Sommerfeld equation $\OS_{\alpha,\lambda} (\phi) =0$ which have either
a ``slow behavior'' $e^{\pm\mu_s z}$ 
or a ``fast behavior'' $e^{\pm \mu_f^+ z}$ as $z \to + \infty$. The two slow modes appear to be perturbations of solutions
of the Rayleigh equation 
$$
Ray_{\alpha,\lambda}(\phi) = (\lambda + i \alpha U)\Delta_\alpha\phi - i\alpha U''\phi =0,
$$
whereas the two fast modes are linked to the Airy type equation
$$(\lambda + i\alpha U- \nu \Delta_\alpha
)\Delta_\alpha \phi =0,$$
or recalling $\mu_f$ introduced in \eqref{def-sfrate}, 
\beq \label{Air}
\nu (\partial_z^2 - \mu_f^2)\Delta_\alpha \phi =0.
\eeq
Let us first consider  the Rayleigh equation $Ray_{\alpha,\lambda}(\phi) = 0$. As $z$ goes to $+\infty$, this equation "converges"
to $\Delta_\alpha \phi = 0$, hence $Ray_{\alpha,\lambda}(\phi) = 0$ has two solutions $\phi_{\alpha,\pm}$, with respective
behaviors $e^{\pm |\alpha| z}$ at infinity. We define the Evans function $E(\alpha,\lambda)$ by
\begin{equation}\label{def-EvansE}
E(\alpha,\lambda) = \phi_{\alpha,-}(0) .
\end{equation}
Note that the Rayleigh equation degenerates at points where $\lambda + i\alpha U (z)$ vanishes. 
In this paper, we restrict ourselves to the case when $\lambda$ is away from the range of $-i\alpha U$. 
Precisely, throughout the paper, letting $\epsilon_0$ be an arbitrarily small, but fixed, positive constant, 
we shall consider the range of $(\alpha,\lambda)$ in $\RR\setminus\{0\}\times \CC$ so that 
\begin{equation}\label{range-c} 
d(\alpha,\lambda) =\inf_{z\in \RR_+}| \lambda + i \alpha U(z)| \ge \epsilon_0.
\end{equation}
Note that $d(\alpha,\lambda) = \Re \lambda$ if $\Im \lambda \in - \alpha \mathrm{Range}(U)$. In any case,
\beq \label{range-c1}
d(\alpha,\lambda) \ge | \Re \lambda | .
\eeq
It turns out that two independent "slow" solutions of Orr Sommerfeld equations can be constructed
as perturbations of these two solutions of Rayleigh equation.

The two "fast" solutions come from the Airy equation \eqref{Air}. This equation degenerates when $\lambda + \alpha^2 \nu + i\alpha U$ gets small.
Points $z_c$ such that $\alpha U(z_c) = - \Im \lambda$ are called "critical layers". The behavior of Airy equation changes
as we approach these points, and in this paper we only study this equation away from these critical layers.
Let us quantify this notion. The Airy's equation has a typical length scale
$$
\delta(z)  = \frac{1}{\mu_f(z)}= \sqrt{\nu \over \lambda + \nu \alpha^2 + i \alpha U(z)} .
$$
If $\delta(z)$ varies within a length $\delta(z)$, namely if $\delta'(z) \delta(z) \sim \delta(z)$, or equivalently 
if $\alpha \sim \nu^{-1/2}$ then the nature of the construction changes (see \eqref{singu} for more details).
In this paper we restrict to the case $| \alpha | \ll \nu^{-1/2}$ or more precisely on $| \alpha | \le \nu^{-\zeta}$
for some $\zeta < 1/2$. 

%

We are mainly interested in getting bounds on the Green function when $\lambda$ has a small positive real part. 
In this case, the condition \eqref{range-c} implies
\begin{equation}\label{comp-musf} \Re \mu_f(z) = \nu^{-1/2} \Re \sqrt{\lambda + \nu \alpha^2 + i\alpha U(z)} \ge \nu^{-1/2} \sqrt{\epsilon_0/2} \gg \mu_s\end{equation}
for sufficiently small $\nu$ and for $| \alpha | \le \nu^{-\zeta}$
for some $\zeta < 1/2$. We may also use these Green function bounds in order to obtain bounds on the solutions of linearized
Navier Stokes equations, through contour integrations. It turns out that these contours may be chosen 
such that $\mu_s \le \Re \mu_f$. Therefore,  we focus on this case in this paper, leaving
aside the case when $\mu_s / \mu_f \ge 1$.

Our main result is the following. 

\begin{theorem}\label{theo-GreenOS}
Let $U(z)$ be a 
boundary layer profile which satisfies \eqref{def-Ubl}. 
For each $\alpha,\lambda$, let by $G_{\alpha,\lambda}(x,z)$ be
the Green kernel of the Orr-Sommerfeld equation, with source term in $x$, and let 
\begin{equation}\label{def-mMf}
\mu_s = | \alpha| , \qquad   \mu_f(z) =  \nu^{-1/2} \sqrt{\lambda + \nu \alpha^2 + i \alpha U(z)}, 
\end{equation}
where we take the square root with positive real part. Let $0 < \theta_0 < 1$ and $\zeta < 1/2$. Let $\sigma_0 > 0$ be arbitrarily small.
Then, there exists $C_0 > 0$ so that 
\begin{equation}\label{est-GrOS}
 |G_{\alpha,\lambda}(x,z)| 
  \le \frac{C_0}{  \mu_s  d(\alpha,\lambda) }e^{-\theta_0 \mu_s |x-z|} 
  +    \frac{C_0}{ |\mu_f(x) |  d(\alpha,\lambda)}  e^{- \theta_0 | \int_x^z \Re \mu_f \; dy|} 
\eeq
uniformly for all $x,z\ge 0$ and $0 < \nu \le 1$, and uniformly in $(\alpha,\lambda)\in \RR\setminus\{0\}\times \CC$ so that $| \alpha | \le\nu^{-\zeta}$, \eqref{range-c} holds, and
$$
| E(\alpha,\lambda)| > \sigma_0 .
$$
In particular, we have
\begin{equation}\label{est-GrOSsimply}
 |G_{\alpha,\lambda}(x,z)| 
  \le \frac{C_0}{  \mu_s  | \Re \lambda |}e^{-\theta_0 \mu_s |x-z|} 
  +    \frac{C_0}{ |\mu_f(x) |  | \Re \lambda |}  e^{- \theta_0 | \int_x^z \Re \mu_f \; dy|} .
\eeq
In addition, there hold the following derivative bounds
\begin{equation}\label{est-GrOS-d}
 | \partial_x^k \partial_z^\ell G_{\alpha,\lambda}(x,z)|  
 \le \frac{C_0\mu_s^{k+\ell}}{   \mu_s  d(\alpha,\lambda)}  e^{-\theta_0 \mu_s |x-z|} 
 +    \frac{C_0 | \mu_f(y) |^{k+\ell}}{| \mu_f(x)| d(\alpha,\lambda) }  e^{- \theta_0 | \int_x^z \Re \mu_f \; dy|}  
\eeq
for all $x,z\ge 0$ and $k, \ell \ge 0$, in which $M_f = \sup_z \Re  \mu_f(z) $. 
Moreover, 
\begin{equation}\label{est-GrOS-delta}
 | \Delta_\alpha G_{\alpha,\lambda}(x,z)|  
 \le \frac{C_0}{   d(\alpha,\lambda) } e^{-\theta_0 \mu_s |x-z|} 
 +    \frac{C_0 | \mu_f(y)|^2}{| \mu_f(x)| d(\alpha,\lambda) }   e^{- \theta_0 | \int_x^z \Re \mu_f \; dy|}  
\eeq
where we "gain" a factor $\mu_s$ in the first term on the right hand side.
\end{theorem}
We believe that the $\theta_0$ factor is purely technical, and that this Theorem holds true for $\theta_0 = 1$. In addition, we note that $\Delta_\alpha G_{\alpha,\lambda}$ enjoys better bounds since $\Delta_\alpha e^{\pm |\alpha| z}Ê= 0$.

To prove this Theorem we first construct approximate solutions to the Orr Sommerfeld equation, and then construct
an approximate Green function. An iteration argument yields the exact Green function together with the stated bounds. 
Our construction of the Green function for the Orr-Sommerfeld problem  was inspired by the pointwise Green function 
approach introduced by Zumbrun-Howard \cite{ZH} and Zumbrun \cite{Z1,Z2}. 

We are also interested in the construction of a pseudo inverse of Orr Sommerfeld operator near a simple eigenvalue, a construction
which is detailed in Section \ref{sectionpseudo}.


\section{Approximate solutions of Orr-Sommerfeld}


In this section, we  construct  four independent approximate solutions to the Orr Sommerfeld equations $\OS_{\alpha,\lambda}(\phi) =0$, 
two with a "fast" behavior and two with a "slow" one. 
The fast modes are constructed using geometrical optics methods, namely following the BKW method.
 For the slow modes we will distinguish between three regimes:

\begin{itemize}

\item bounded $| \alpha |$. In this case the  slow modes are perturbations of the eigenmodes of Rayleigh equations.

\item $ 1 \ll | \alpha | \le \nu^{-1/4}$ (or any small negative power of $\nu$). 
We use the fact that Rayleigh equation is a perturbation of $\Delta_\alpha$.
The slow modes are perturbations of the eigenmodes of $\partial^2_z - \alpha^2$, namely
$e^{\pm |\alpha| z}$.

\item $ \nu^{-1/4} \le | \alpha | \le \nu^{-\zeta}$, for $\zeta <1/2$.  
In this case $e^{\pm |\alpha| z}$ is a sufficient approximation.

\end{itemize}
Solutions will be constructed in function spaces $L^\infty_\eta$, for $\eta>0$, that consist of smooth functions $f$ so that the norm 
$$
\|f\|_\eta: = \sup_{z\ge 0} e^{\eta |z|} |f(z)| 
$$ 
is finite.


\subsection{Fast modes}\label{sec-fast}


In this section, we shall construct two independent approximate solutions, which asymptotically behave like $e^{\pm \mu^+_f z}$, 
of the Orr-Sommerfeld equation $\OS_{\alpha,\lambda}(\phi) = 0$.  We will use the BKW method. Let us first discuss its validity.
Note that locally the characteristic length scale of the oscillations is
$$
\delta(z) = {1 \over \mu_f(z)} =\sqrt{\nu \over \lambda + \nu \alpha^2 + i \alpha U(z) }.
$$
The BKW method is valid provided $\delta$ has a small change during a period, namely provided
$
\delta' \delta \ll \delta,
$
or equivalently
\beq \label{singu}
\delta'(z) = {-i\sqrt \nu \alpha U'(z) \over 2 (\lambda + \nu \alpha^2 + i \alpha U(z))^{3/2} } \ll 1.
\eeq
Note that it may happen that for some particular $z_c$, $\Im \lambda + \alpha U(z_c) = 0$. Such $z_c$ are called
critical layers, or turning points.
If $\alpha \sim \nu^{-1/2}$, then the denominator and numerator of (\ref{singu}) are of order $O(1)$ at such points, hence
the condition (\ref{singu}) is not satisfied and $\delta'(z_c) \sim 1$. 

On the contrary if $\alpha \lesssim \nu^{-\zeta}$ with $\zeta < 1/2$, then near critical points, the denominator is of
order $O(1)$ but the numerator is of order $O(\nu^{1/2 - \zeta})$. Therefore the condition (\ref{singu}) is satisfied
provided $\nu$ is small enough. Similarly, $\mu_f^{-1-j}\partial_z^j \mu_f(z)$ is of order $O(\nu^{1/2 - \zeta})$ or smaller for $j \ge 1$.


%
%

\begin{prop}\label{prop-fastOrrapp} 
Let $N > 0$ be arbitrarily large. Then for sufficiently small $\nu$ and for $|\alpha|\le \nu^{-\zeta}$ with $\zeta <1/2$, 
there exist two approximate  solutions $\phi_{f,\pm}^{app}(z)$ which solve Orr-Sommerfeld equations up to a small error term
$$
\OS_{\alpha,\lambda}(\phi_{f,\pm}^{app}) = O(\nu^N | \phi_{f,\pm}^{app}|),
$$
with $\phi_{f,\pm}^{app}(0) = 1$ and
\begin{equation}\label{fast-mode1} 
\phi_{f,\pm}^{app} (z) =  e^{\pm \int_0^z \mu_{f}(y) \; dy  } \Big( 1 + \phi_{\pm}(z)\Big) ,
\end{equation}
where $\phi_{\pm}$ and their derivatives are uniformly bounded in $\alpha$, $\nu$ and $z$, 
and converge exponentially fast to $0$ at $z=+ \infty$.
\end{prop}

\begin{proof}
Following a semi classical approach, we look for $\phi_{f,\pm}^{app}$ under the form
$$
\phi_{f,\pm}^{app} = \exp \Bigl({\theta_{\pm}^{app} \over \sqrt{\nu}} \Bigr) .
$$
Let $\theta = \theta_{\pm}^{app}$ to simplify the notations. We compute 
$$
\partial_z^2 \phi_{f,\pm}^{app} = \Bigl( {\theta'^2 \over \nu} + {\theta'' \over \sqrt{\nu}} \Bigr) \phi_{f,\pm}^{app} 
$$
and
$$
\nu \partial_z^4 \phi_{f,\pm}^{app} 
= \Bigl( {\theta'^4 \over \nu} + 6 {\theta'^2 \theta'' \over \sqrt\nu} 
+ 4 \theta' \theta''' + 3 \theta''^2 + \sqrt\nu \theta'''' \Bigr) \phi_{f,\pm}^{app} .
$$
We now expand $\theta$ in powers of $\sqrt\nu$; namely, 
$$
\theta = \sum_{i = 0}^N \theta_j \nu^{j/2},
$$
where the functions $\theta_j$ will themselves depend on $\alpha$ and $\lambda$. Putting the Ansatz into the Orr-Sommerfeld equations, at leading order, we obtain
$$
(\lambda + i\alpha U ) (\theta_0'^2 - \nu \alpha^2) - \Bigl( \theta_0'^4 - 2  \nu \alpha^2 \theta_0'^2 + \nu^2 \alpha^4  \Bigr) = 0.
$$
Factorizing by $\theta_0'^2 - \nu \alpha^2$ we get
$$
\theta_0'^2 =  \lambda + \nu \alpha^2 + i \alpha U = \nu \mu_f^2(z),
$$
which gives
$$
\theta_0' = \pm \sqrt{\nu}  \mu_f(z).
$$
Note that $\theta_0'$ converges exponentially fast to $\pm \sqrt{\nu} \mu_f^+$ and $\theta_0''$ converges exponentially fast
to $0$. To obtain $\theta_1$ we equate the powers in $\sqrt{\nu}^{-1}$ and get
$$
- 4  \theta_0'^3 \theta_1' + 4  \nu \alpha^2 \theta_0' \theta_1' + 2 ( \lambda + i\alpha U ) \theta_0' \theta_1' = S,
$$ 
where the source term $S = 6 \theta_0'^2 \theta_0''$ only depends on $\theta_0'$ and its derivatives.
This leads to
$$
\theta_1' =   {S  \over (- 4 \theta_0'^2  +4  \nu \alpha^2  + 2 (\lambda + i\alpha U ) ) \theta_0'}
= -{S  \over 2 (\lambda + i\alpha U ) \theta_0'} .
$$
As $\theta_0'$ is bounded away from $0$, $\theta_1'$ is correctly defined. Moreover, 
$\theta_1$ converges exponentially at infinity, as well as all its derivatives, and as 
$\theta_0'' = O(\alpha)$, $\theta_1 = O(\alpha)$. This leads to
\beq \label{bornealpha}
\theta^{app}_\pm = \theta_0 + O(\alpha \nu^{-1/2}).
\eeq
We then obtain equations and similar estimates on the remaining $\theta_j$ by equaling successive powers of $\nu$.
 The Proposition follows.
\end{proof}


\subsection{Slow modes}


\begin{prop} \label{slow1}
There exist
two solutions $\phi_{s,\pm}^{app}$ which approximately solve the Orr Sommerfeld equations: precisely, for any $N$, 
$$
| \OS_{\alpha,\lambda}(\phi_{s,\pm}^{app}) | \le C_N \nu^N e^{\pm | \alpha | z - \eta z} 
$$
and behave like $e^{\pm | \alpha | z}$ as $z$ goes to $+ \infty$: for any $n$, 
$$
| \partial_z^n \phi_{s,\pm}^{app} (z) | \le C_n e^{\pm | \alpha | z} . 
$$
\end{prop}
For the proof of Proposition \ref{slow1}, we shall distinguish three cases: bounded $\alpha$, moderate $\alpha$, 
and large $\alpha$, that will be detailed in the next
sections. We restrict ourselves to $\alpha > 0$, the opposite case being similar.


\subsubsection{Approximate slow modes for bounded $\alpha$ and $\lambda$}


As $z$ goes to $+ \infty$, the Rayleigh equation "converges" to $\Delta_\alpha \phi$. 
Therefore the Rayleigh equation admits two particular equations, called $\phi_{\alpha,\pm}$ which behave like
$e^{\pm | \alpha | z}$ as $z \to + \infty$. Moreover $|\partial^n_z \phi_{\alpha,\pm}(z) | \le C_n e^{\pm | \alpha | z}$ for every
positive $n$. Note that
$$
\OS_{\alpha,\lambda}(\phi_{\alpha,\pm}) = - \nu \Delta_\alpha^2 \phi_{\alpha,\pm} .
$$
Using the Rayleigh equation, we compute 
$$
\Delta_\alpha \phi_{\alpha,\pm} = {i\alpha U'' \phi_{\alpha,\pm} \over \lambda + i\alpha U},
$$
which gives
$$
\OS_{\alpha,\lambda}(\phi_{\alpha,\pm}) =  - \nu \Delta_\alpha \Bigl( {i\alpha U'' \phi_{\alpha,\pm} \over \lambda + i\alpha U} \Bigr) 
$$
$$
 = - \nu  \Bigl( {i\alpha U''  \over \lambda + i\alpha U} \Bigr)^2 \phi_{\alpha,\pm}
 - 2  \nu  \partial_z \phi_{\alpha,\pm} \partial_z \Bigl( {i\alpha U'' \over \lambda + i\alpha U} \Bigr) 
- \nu \phi_{\alpha,\pm} \partial_z^2  \Bigl( {i\alpha U'' \over \lambda + i\alpha U} \Bigr) .
$$
Note that $\lambda + i\alpha U$ is bounded away from $0$, therefore 
$$
| \OS_{\alpha,\lambda}(\phi_{\alpha,\pm}) | \le C \nu e^{\pm | \alpha | z - \eta z} ,
$$
and similarly for all its derivatives.

We now look for approximate solutions of Orr Sommerfeld solutions $\phi_{s,\pm}^{app}$ of the form
$$
\phi_{s,\pm}^{app} = \sum_{j=0}^N \phi_{\alpha,\pm}^j
$$
for arbitrarily large $N$, starting with $\phi_{\alpha,\pm}^0 = \phi_{\alpha,\pm}$.
We have
$$
Ray_\alpha(\phi_{\alpha,\pm}^{j+1}) = - \OS_{\alpha,\lambda}(\phi_{\alpha,\pm}^j) .
$$
Note that 
\begin{equation}
\OS_{\alpha,\lambda}(\phi_{s,\pm}^{app}) = - \nu \Delta_\alpha^2 \phi_{\alpha,\pm}^N .
\end{equation}
We will focus on the construction of $\phi_{s,-}^{app}$, the construction of $\phi_{s,+}^{app}$ being similar.
To end the proof of Proposition \ref{slow1} we need to bound the various $\phi_{\alpha,-}^i$,
which is done through the iterative use of the following Proposition.

\begin{prop} 
There exist constants $C_n$ such that the following assertion is true.
For any $\beta>0$ and any smooth function $\psi$ there exists a smooth solution $\phi$ of $Ray_\alpha(\phi) = \psi$
such that 
$$
\sup_{k \le n} \| \partial_z^k \phi \|_{\alpha} 
+ \sup_{k \le n} \| \partial_z^n \Delta_\alpha \phi \|_{\alpha + \beta}
\le {C_n \over E(\alpha,\lambda)} \sup_{k \le n} \| \partial_z^k \psi \|_{\alpha+\beta}
$$
where $\|\phi\|_\eta = \sup_{z\ge 0} e^{\eta |z|}|\phi(z)|$. \end{prop}

\begin{proof}
We first construct the Green function of the Rayleigh operator. Let
$$
\widetilde \phi_{\alpha,+}(z) = \phi_{\alpha,-}(0) \phi_{\alpha,+}(z)
- \phi_{\alpha,+}(0) \phi_{\alpha,-}(z)
$$
Then $\widetilde \phi_{\alpha,+}(0) = 0$ and the Wronskian of $\widetilde \phi_{\alpha,+}$ and 
$\phi_{\alpha,-}$ equals
$$
W(\widetilde \phi_{\alpha,+}, \phi_{\alpha,-}) = \phi_{\alpha,-}(0) W(\phi_{\alpha,+},\phi_{\alpha,-})
= 2 \alpha \phi_{\alpha,-}(0)
$$
evaluating this latest Wronskian at infinity. The Green function of the Rayleigh operator is therefore
$$
G(x,z) = {1 \over 2 \alpha \phi_{\alpha,-}(0)} \phi_{\alpha,-}(x) \widetilde \phi_{\alpha,+}(z) \qquad
\hbox{if} \qquad z < x
$$
$$
G(x,z) = {1 \over 2 \alpha \phi_{\alpha,-}(0)} \widetilde \phi_{\alpha,+}(x) \phi_{\alpha,-}(z) \qquad
\hbox{if} \qquad z  > x.
$$
We then have
$$
\phi(z) = \int_0^{+ \infty} G(x,z) \psi(x) dx .
$$
Using the asymptotic behavior of $\phi_{\alpha,\pm}$ we get the claimed bounds on $\| \partial_z^n \phi \|_{\alpha}$
with $n = 0$ and $n = 1$ by a direct computation. Higher derivatives are obtained by differentiating
$$
\partial_y^2 \phi = \alpha^2 \phi + {i\alpha U'' \over \lambda + i\alpha U} \phi + \psi,
$$
keeping in mind that $\alpha$ is bounded and $\lambda$ is away from the range of $-i\alpha U$. 
Next, we write 
$$
\Delta_\alpha \phi = {i\alpha U'' \over \lambda + i\alpha U} \phi + \psi
$$
which gives the desired bounds on $\Delta_\alpha \phi$.
\end{proof}


\subsubsection{Approximate slow modes for $1 \ll |\alpha| \le \nu^{-1/4}$ or large $\lambda / \alpha$}


For large $\alpha$, or for large $\lambda / \alpha$,
 the Rayleigh operator is a small perturbation of $\partial_z^2 - \alpha^2$
 and we can construct approximate eigenmodes $\phi_{s,\pm}^{app}$ using a perturbative construction.
Namely, the Rayleigh equation may be rewritten as 
$$ 
\Delta_\alpha \phi = { i\alpha U'' \phi \over  \lambda + i \alpha U} .
$$
Note that $\alpha^{-1} e^{-\alpha | x - z|}$ is a Green function for $\Delta_\alpha$.
We therefore define the following operator $\mathcal{T}$ by
$$ 
\mathcal{T}[\phi](z) := \int_0^\infty \alpha^{-1} e^{-\alpha |x-z|}  {i \alpha U'' \phi(x) \over \lambda + i \alpha U } \; dx .
$$
We shall prove that for sufficiently large $\alpha$, the map $\mathcal{T}$ is well-defined and contractive 
from $L^\infty_{\alpha +\eta}$ to itself. 
Indeed, for $\phi \in L^\infty_{\alpha + \eta}$, as $\lambda + i\alpha U$ is bounded away from $0$, we have 
$$
 | \mathcal{T}[\phi](z)| 
 \le C_0 \int_0^\infty  e^{-\alpha |x-z|} e^{-\eta x - \alpha x}   \| \phi\|_{\alpha + \eta} \; dx 
\le C_0 \alpha^{-1} \| \phi\|_{\alpha + \eta}  e^{-\eta z - \alpha z}. 
$$
This proves that $\mathcal{T}[\phi] \in L^\infty_{\alpha + \eta}$. 
If $\alpha$ is large enough then $\mathcal{T}$ is a contraction in this space. On the other hand, if $\lambda / \alpha$ is large enough we rewrite
$$
 {i\alpha U'' \phi(x) \over  \lambda + i\alpha U} =  { U'' \phi(x) \over  U- i \alpha^{-1} \lambda }
 $$
 which is bounded by $C / (\alpha^{-1} \lambda)$. Hence $\mathcal{T}$ is a contraction if $\lambda / \alpha$ is large enough.
 
We now construct two independent solutions of the Rayleigh equation, which behaves like
$e^{\pm \alpha z}$ for large $z$. Let us detail the "-" case. We look for $\phi_{s,-}$ under the form
$$
\phi_{s,-} = \sum_{n \ge 0} \phi_{-}^n
$$
with $\phi_-^0 = e^{- \alpha z}$ and 
$\phi_{-}^{n+1} = \mathcal{T}[\phi_{-}^n]$. As $\mathcal{T}$ is contractive, the previous sum converges in $L^{\infty}_{\alpha + \eta}$.
Note that in particular
$$
\phi_{\alpha,-} = e^{- \alpha z}  ( 1 + O(\alpha^{-1})_{L^{\infty}_{\alpha+\eta}}),
$$ 
and similarly for its derivatives. The construction of $\phi_{\alpha,+}$ is similar.

The construction of approximate solutions of Orr Sommerfeld is similar to that of the previous section. We start with
$\phi_{s,-}$ and note that
$$
 \nu \| \Delta_\alpha^2 \phi_{s,-} \|_{\alpha + \eta} \le C \nu | \alpha |^2 \lesssim \nu^{1/2}.
$$
We then introduce $\phi_{s,-}^1$, defined by
$$
Ray_\alpha(\phi_{s,-}^1) = - \nu \Delta_\alpha^2 \phi_{s,-} ,
$$
which can be bounded using the $\mathcal{T}$ operator.  To  end the proof of Proposition \ref{slow1},
we iterate the construction as in the previous section.


\subsubsection{Approximate slow modes for $\nu^{-1/4} \le |\alpha|  \ll \nu^{-1/2}$}


We look for eigenmodes of the form
$$
\phi_{s,\pm}^{app} = \exp ( \alpha \theta_{\pm}^{app} )
$$
where $\theta_\pm^{app}$ may be expanded in powers of $\alpha^{-1}$.
As in Section \ref{sec-fast}, we get
$$
- \nu \alpha^4 \theta_0'^4 +2  \nu \alpha^4 \theta_0'^2 - \nu \alpha^4  
+ (\lambda + i\alpha U) (\alpha^2 \theta_0'^2 - \alpha^2) = 0,
$$
This time we choose $\theta_0 = \pm 1$ and iterate as in Section \ref{sec-fast} to prove Proposition \ref{slow1}.
Note again that the leading order of $\Delta_\alpha \phi_{s,\pm}^{app}$ vanishes.


\section{Approximate Green function}\label{sec-Greenapp}


We now construct an approximate Green function $H^{app}$ using the approximate solutions
$\phi_{s,\pm}^{app}$ and $\phi_{f,\pm}^{app}$.
We will decompose this Green function into two components
$$
H^{app} = G^{app} + \hat G^{app}
$$
where $G^{app}$ does not take into account the boundary conditions and focus on the discontinuity at $y = x$,
and where $\hat G^{app}$ restores the proper boundary conditions.

Hence, first forgetting the boundary condition, we look for $G^{app}(x,y)$ of the form
\begin{equation}\label{def-GappX}
\begin{aligned}
G^{app}(x,y) &= a_+(x)  {\phi_{s,+}^{app}(y) \over c_2}
+ {b_+(x) }  {\phi_{f,+}^{app}(y) \over \phi_{f,+}^{app}(x)} 
\quad \hbox{for} \quad y < x,
\\
G^{app}(x,y) &= a_-(x)  {\phi_{s,-}^{app}(y) \over c_1} 
+ {b_-(x)} {\phi_{f,-}^{app}(y) \over \phi_{f,-}^{app}(x)} 
\quad \hbox{for} \quad y  > x,
\end{aligned}\end{equation}
where the normalization constants $c_1$ and $c_2$ will be fixed later.
Let 
\begin{equation}\label{def-vvvapp}
v(x) = (- a_-(x), a_+(x), - b_-(x), b_+(x) ) .
\end{equation}
By definition, $G^{app}$, $\partial_y G^{app}$, $\sqrt{\nu} \partial_y^2 G^{app}$ are continuous at $x = y$
and $\nu \partial_y^3 G^{app}$ has a jump at $x = y$, of magnitude $1$. 
Let
\beq \label{matriceM}
M = \left( \begin{array}{cccc} 
\phi_{s,-} / c_1& \phi_{s,+} / c_2 & \phi_{f,-}  & \phi_{f,+}  \cr
\partial_y \phi_{s,-} / c_1 \mu_f & \partial_y\phi_{s,+} / c_2 \mu_f
& \partial_y\phi_{f,-} /  \mu_f & \partial_y\phi_{f,+} / \mu_f  \cr
\partial_y^2 \phi_{s,-} / c_1 \mu_f^2 &\partial_y^2 \phi_{s,+} / c_2 \mu_f^2 
& \partial_y^2 \phi_{f,-} /  \mu_f^2  & \partial_y^2 \phi_{f,+} /  \mu_f^2 \cr
 \partial_y^3 \phi_{s,-} / c_1 \mu_f^3 &  \partial_y^3 \phi_{s,+} / c_2  \mu_f^3
& \partial_y^3 \phi_{f,-}  /  \mu_f^3 &  \partial_y^3 \phi_{f,+} /  \mu_f^3 \cr 
\end{array} \right) ,
\eeq
where the functions $\phi_{s,\pm} = \phi_{s,\pm}^{app}$ and $\phi_{f,\pm} = \phi_{f,\pm}^{app}$ and their derivatives are evaluated at $y=x$.
Then 
\beq \label{Mv}
M v = (0,0,0,1/ \nu \mu_f^3) .
\eeq
In the following sections we will bound the solution $v$ of (\ref{Mv}).
Let us define the four two by two matrices $A$, $B$, $C$ and $D$  by
$$
M = \left( \begin{array}{cc} 
A & B \cr
C & D \cr \end{array} \right) .
$$
Note that, using (\ref{bornealpha}),
$$
D = \left( \begin{array}{cc}
1  & 1 \cr
- 1 & 1 \cr 
\end{array} \right) + O(\alpha \mu_f^{-1}).
$$
 Hence the matrix $D$ is bounded and invertible, upon recalling that $\alpha \ll \mu_f$ in the range of $\alpha$ that we consider
 (see \eqref{comp-musf}). Moreover its inverse is bounded and equals
$$
D^{-1} = {1 \over 2}  \left( \begin{array}{cc}
1  & - 1 \cr 
1 & 1 \cr 
\end{array} \right) + O(\alpha \mu_f^{-1}).
$$
We shall consider two cases: bounded $\alpha$ and unbounded $\alpha$. 


\subsection{First case: 
bounded $\alpha$}\label{sec-bda}


We take $c_1 = c_2 = 1$. Note that $A = A_1 A_2$
where
$$
A_1 = \left( \begin{array}{cc}  
1 & 0 \cr
0 & \mu_f^{-1} \cr
\end{array} \right),
\quad 
A_2 = \left( \begin{array}{cc}  
\phi_{s,-} & \phi_{s,+} \cr
\partial_y \phi_{s,-} & \partial_y \phi_{s,+} \cr
\end{array} \right) .
$$
The determinant $E^{app}(\alpha,\lambda)$ of $A_2$ is a perturbation of the Evans function $E(\alpha,\lambda)$
in the sense
$$
E^{app}(\alpha,\lambda) = E(\alpha, \lambda) + O(\nu^\sigma),
$$
for some positive $\sigma$.
Hence if $E(\alpha, \lambda) \ne 0$, then $A_2$ and $A$ 
are invertible provided $\nu$ is small enough, and $A_2^{-1}$ is bounded. 
Moreover the matrix $M$ has an approximate inverse
$$
\widetilde M = \left( \begin{array}{cc}
A^{-1} & - A^{-1} B D^{-1}  \cr
0 & D^{-1} \cr 
\end{array} \right) 
$$
in the sense that $M  \widetilde M = Id + N$ where
$$
N =  \left( \begin{array}{cc}
0 & 0 \cr 
 C A^{-1} &  - C A^{-1} B D^{-1} \cr 
\end{array} \right) .
$$
Note that $C$ is of order $O(\mu_f^{-2})$ since $\alpha$ is bounded, that $B$ is bounded 
and that $A^{-1} = A_2^{-1} A_1^{-1}$ is of order $O(\mu_f)$.
Hence we have  $N = O(\mu_f^{-1})$. Therefore $(Id + N)^{-1}$
is well defined and uniformly bounded for $\nu$ small enough provided $E(\alpha,\lambda) \ne 0$.
As a consequence, 
$$
M^{-1} = \widetilde M (Id + N)^{-1} = \widetilde M \sum_n N^n.
$$
Note that the two first lines of $N^n$ vanish. Therefore
$$
(Id + N)^{-1} (0,0,0, 1 / \nu \mu_f^3) = \Bigl( 0, 0, O(1 / \nu \mu_f^4), 1 / \nu \mu_f^3  + O(1 / \nu \mu_f^4) \Bigr) .
$$
As $D^{-1}$ is bounded and $A^{-1} B D^{-1}$ is of order $O(\mu_f)$, we obtain that 
$a_\pm$ and $b_\pm$ are respectively of order $O(1 / \nu \mu_f^2)$ and $O(1/ \nu \mu_f^3)$.
Note that $\alpha$ is bounded in this case, which give the desired bounds since
$$
\nu \mu_f^2 = \lambda + \nu \alpha^2 + i  \alpha U 
$$
hence
$$
| \nu \mu_f^2 | \ge d(\alpha,\lambda),
$$
which ends this first case.


\subsection{Case $2$: 
large $\alpha$}


We take $c_1 = \phi_{s,+}^{app}(x)$ and $c_2 = \phi_{s,-}^{app}(x)$.
In this case $A$ is of the form
$$
A = \left( \begin{array}{cc}
1 & 1 \cr
- \alpha \mu_f^{-1} & \alpha \mu_f^{-1} \cr \end{array} \right) (1 + o(1)).
$$
Its inverse $A^{-1}$ equals
$$
A^{-1} = {1 \over 2} \left( \begin{array}{cc}
1 & - \alpha^{-1} \mu_f \cr
1 &  \alpha^{-1} \mu_f \cr \end{array} \right) (1 + o(1)).
$$
Note that $D^{-1}$ and $B$ are bounded and $A^{-1}$ is order $O(\mu_f / \alpha)$. 
As $C$ is of order $O(\alpha^2 / \mu_f^{2})$, $N$ (defined in the previous section) is of order $O(\alpha / \mu_f)$.
Hence, as $| \alpha |\ll \mu_f$ in view of \eqref{comp-musf}, we have 
$$
(Id + N)^{-1} = \sum_n (-1)^n N^n.
$$
This leads to
\beq \label{IdN1}
(Id + N)^{-1}(0,0,0,1/\nu \mu_f^4) =  \Bigl( 0,0,O(\alpha / \nu \mu_f^4), O(1 / \nu \mu_f^3) \Bigr) .
\eeq
It remains to evaluate the image of this vector by $\widetilde M$. As $D^{-1}$ is bounded, we obtain that
$b_{\pm}$ are of order $O(1 / \nu \mu_f^3) = O(1 / \mu_f d(\alpha,\lambda))$.

Moreover, we compute 
$$
D^{-1}(0, O(1 / \nu \mu_f^3)) = \Bigl[ (-1,1) + O(\alpha \mu_f^{-1}) \Bigr]ÊO(1 / \nu \mu_f^3).
$$
As
$$
B = \left( \begin{array}{cc}
1 & 1 \cr
-1 & 1 \cr
\end{array} \right) (1 + O(\mu_f^{-1})),
$$
we obtain
$$
B D^{-1}(0, O(1 / \nu \mu_f^3)) =  \Bigl[ (0,1) + O(\alpha \mu_f^{-1}) \Bigr]ÊO(1 / \nu \mu_f^3).
$$
As a consequence, we obtain 
$$
A^{-1} B D^{-1} (0, O(1 / \nu \mu_f^3)) = O(1 / \alpha \nu \mu_f^2).
$$
It remains to bound the images of the $O(\alpha / \nu \mu_f^4)$ term in the equation (\ref{IdN1}). We have 
$$
D^{-1}(O(\alpha / \nu \mu_f^4),0) = \Bigl[ (1,1) + O(\alpha \mu_f^{-1}) \Bigr]ÊO(\alpha / \nu \mu_f^4).
$$
Hence
$$
B D^{-1} (O(\alpha / \nu \mu_f^4),0)  =  \Bigl[ (1,0) + O(\alpha \mu_f^{-1}) \Bigr]ÊO(\alpha / \nu \mu_f^4)
$$
and $A^{-1} B D^{-1} (O(\alpha / \nu \mu_f^4),0) = O(\alpha / \nu \mu_f^4)$. Using again $\alpha \ll \mu_f$, we obtain
that $a_{\pm}$ are of order $O(1 / \nu \mu_f^2 \alpha) = O(1 / \alpha d(\alpha,\lambda))$.


\subsection{Boundary condition}


We now add to $G^{app}$ another approximate Green function $\hat G^{app}$ to handle the boundary conditions.
We look for $\hat G^{app}$ under the form
$$
\hat G^{app}(y) = d_s {\phi_{s,-}(y) \over d_1} + d_f {\phi_{f,-}(y)  \over \phi_{f,-}(0)}
$$
where the normalization constant $d_1$ will be fixed later, 
and look for $d_s$ and $d_f$ such that
\beq \label{Greenb1}
G^{app}(x,0) + \hat G^{app}(0) = 
\partial_y G^{app}(x,0) + \partial_y \hat G^{app}(0) =  0.
\eeq
Let
$$
\hat M = \left( \begin{array}{cc} \phi_{s,-} / d_1 & \phi_{f,-} / \phi_{f,-}(0) \cr
\partial_y \phi_{s,-} / d_1 & \partial_y \phi_{f,-} / \phi_{f,-}(0) \cr \end{array} \right) ,
$$
the functions being evaluated at $y = 0$. Then (\ref{Greenb1}) can be rewritten as
$$
\hat M d = - (G^{app}(x,0), \partial_y G^{app}(x,0)) 
$$
where $d = (d_s,d_f)$. Note that
$$
(G^{app}(x,0), \partial_y G^{app}(x,0))  = Q (a_+,b_+)
$$
where
$$
Q = \left( \begin{array}{cc}
 \phi_{s,+}(0) /c_2  & 1 \cr
 \partial_y \phi_{s,+}(0) / c_2 & \partial_y \phi_{f,+}(0) / \phi_{f,+}(0) \cr
\end{array} \right) .
$$
By construction
\beq \label{defid}
d = -  \hat M^{-1}  Q (a_+,b_+) .
\eeq
Let us first consider bounded $\alpha$. We take $d_1 = 1$. This leads to
$$
\hat M = \left( \begin{array}{cc} \phi_{s,-}(0) & 1 \cr
\partial_y \phi_{s,-}(0) & - \mu_f + O(1) \cr \end{array} \right) .
$$
Note that  $\hat M = M_1 M_2$ with
$$
M_1 = \left( \begin{array}{cc} 1 & 0 \cr 0 &  \mu_f \cr \end{array} \right), \qquad
M_2 = \left( \begin{array}{cc} \phi_{s,-}(0) & 1 \cr
\partial_y \phi_{s,-}(0) / \mu_f & -1 + O( 1 / \mu_f) \end{array} \right) .
$$
The determinant of $M_2$ equals to $-E(\alpha,\lambda) = -\phi_{s,-}(0)$, up to a small term of order $\mu_f^{-1}\sim \sqrt \nu$, recalling that $\alpha$ is bounded.
Hence $M_2$ is invertible, and $M_2^{-1}$ is bounded if $E(\alpha,\lambda) \ne 0$, provided $ \nu$ is small enough.
Then
$$
\hat M^{-1} Q = M_2^{-1} M_1^{-1} Q.
$$
Note that $(a_+,b_+) = (O(1/\nu \mu_f^2), O(1 / \nu \mu_f^3))$. Hence $Q (a_+,b_+) = O(1 / \nu \mu_f^2)$.
Therefore $M_1^{-1} Q (a_+,b_+) = (O(1/ \nu \mu_f^2), O(1/ \nu \mu_f^3))$. Hence, as the second term of the first
column of $M_2$ is of order $O(1 / \mu_f)$ we get, as desired, that 
\beq \label{boundofd}
d = (O(1 / \alpha \nu \mu_f^2), O(1 / \nu \mu_f^3)),
\eeq
keeping in mind that $\alpha$ is bounded.

For large $\alpha$ we choose $d_1 = \phi_{s,-}(0)$. Then
$$
Q = \left(\begin{array}{cc} 1 & 1 \cr
\alpha + O(1) & \mu_f + O(1) \cr
\end{array} \right),
$$
$$
\hat M = \left( \begin{array}{cc} 1 & 1  \cr  - \alpha + O(1)  & - \mu_f + O(1) \cr \end{array} \right)
$$
and
$$
\hat M^{-1} = {1 \over \mu_f - \alpha + O(1)} \left( \begin{array}{cc} \mu_f  + O(1) & 1 \cr
- \alpha + O(1) & -1 \cr \end{array} \right) .
$$
In this case $(a_+,b_+) = (O(1 / \alpha  \nu \mu_f^2),O(1/ \nu \mu_f^3))$. A direct computation of 
$\hat M^{-1} Q (a_+,b_+)$ again gives (\ref{boundofd}).
Combining all the previous estimates ends the proof.


\section{Exact Green function}\label{sec-exactGreen}


Let 
$$
H^{app} = G^{app} + \hat G^{app}
$$ 
be the complete approximate Green function. By construction, $H^{app}$ satisfies the zero boundary conditions \eqref{G-noslip}. 
We now construct the exact Green function $G(x,z)$ as an infinite sum
\beq \label{defG}
G(x,z) = \sum_{n \ge 0} G_n(x,z),
\eeq
where $G_0 = H^{app}$, 
$$
G_1 = - H^{app} \star (\OS_{\alpha,\lambda}(H^{app}) - \delta_{y=x}),
$$ 
and $G_n$ is defined by iteration through
$$
G_{n+1} = - H^{app} \star \OS_{\alpha,\lambda}(G_n).
$$
Hence, it suffices to prove that the series \eqref{defG} converges in a suitable function space, which follows immediately from the following lemma. The stated bounds for $G(x,z)$ in Theorem \ref{theo-GreenOS} then follow from those on $H^{app}(x,z)$. 

\begin{lem}
For each $x$, assume that
$$
| f^x(y) | \le e^{- \alpha'  |x-y | } 
$$
for some $\alpha'$ such that $\alpha' < | \alpha |$ and $\alpha' < \Re \mu_f$.
Then
$$
| \OS_{\alpha,\lambda}(G^{app} \star f^x) (y)| \le C \nu^{N-2} e^{- \alpha' |x-y| } .
$$
\end{lem}
\begin{proof}
Note that
$$
\OS_{\alpha,\lambda}(G^{app} \star f^x)(y)  = \int \OS_{\alpha,\lambda}(G^{app})(z,y) f^x(z) dz .
$$
However we recall that $\phi_{s,\pm}^{app}$ satisfy
$$
| \OS_{\alpha,\lambda}(\phi_{s,\pm}^{app})  | \le C \nu^N e^{\pm |\alpha|  z },
$$
$$
| \OS_{\alpha,\lambda}(\phi_{f,\pm}^{app}) |Ê\le C \nu^N |Ê\phi_{f,\pm}^{app} | .
$$
Using the bounds on the coefficients on $G^{app}(z,y)$, this leads to
$$
| \OS_{\alpha,\lambda}(G^{app}(z,y)) | \le C \nu^{N-2}  e^{- \alpha | y - z |Ê }.
$$
The Lemma follows by convolution.
\end{proof}


\section{Construction of a pseudo inverse \label{sectionpseudo}}


We now focus on the case when $\lambda$ is close to a simple eigenvalue $\lambda_0$.

\begin{theo} \label{theopseudo}
Let $\alpha$ be fixed.
Let $\lambda_0$ be a simple eigenvalue of $Orr_{\alpha,\lambda}$ with corresponding eigenmode
$\phi_{\alpha,\lambda_0}$. Then there exists a bounded family of linear forms $l^\nu$
and a family of pseudoinverse operators $Orr^{-1}_{\alpha,\lambda}$ such that for any stream function $\phi$,
$$
Orr_{\alpha,\lambda}\Bigl(Orr^{-1}_{\alpha,\lambda}(\phi) \Bigr)  = \phi - l^\nu(\phi) \phi_{\alpha,\lambda_0} 
$$
for $\lambda$ near $\lambda_0$. Moreover, the pseudoinverse $Orr_{\alpha,\lambda}^{-1}$ may be defined through a Green function $\widetilde G_{\alpha,\lambda}(x,z)$ which satisfies the same bounds in (\ref{est-GrOS}).
\end{theo}


\subsection{Principle of the construction}


Let us sketch the principle of the proof on a simplified case. Let $A_0$ be a $N \times N$ matrice of rank $N-1$
(which is a toy model for the Rayleigh operator when $\lambda$ is a simple eigenvalue), and let $A(\eps)$
be a bounded family of $N \times N$ matrices (toy model for Orr Sommerfeld equation). We want to construct an inverse for 
$$
A^\eps = A_0 + \eps A(\eps).
$$
Let us first invert $A_0$. Let $v$ be a unit vertor, orthogonal to the image of $A_0$. Let $P$ be the orthogonal projector
on the image of $A$, namely
$$
P v =  f - (f. v) v.
$$
Let $B$ be a pseudo inverse of $A_0$, namely a matrix such that, on the image of $A_0$, $A_0 B =Id$. 
Then  $u = B P f$ solves 
$$
A_0 u = f - (f.v) v .
$$
We now fulfill a similar construction for $A^\eps$ for small $\eps$. 
Let  $u_0 = B  P f$. Then 
$$
A^\eps u_0 = f - (f.v) v  + \eps A(\eps) u_0.
$$
We know define $u_1 = -  B P A(\eps) u_0$. Then $u_0 + u_1$ solves
$$
A^\eps (u_0 + \eps u_1) = f - (f.u_0) v + \eps (A(\eps) u_0. v) v - \eps^2 A(\eps) B P A(\eps) u_0 
$$
and the construction follows by iteration.


\subsection{Rayleigh equation}


In this section we fix $\alpha$ and investigate the Rayleigh operator $Ray_{\alpha,\lambda}$ 
when $\lambda$ is near a simple eigenvalue $\lambda_0$ of $Ray_{\alpha,\lambda}$. 
We will also assume that $Ker(Ray_{\alpha,\lambda_0}^2) = \cit \phi_{\alpha,\lambda_0,\pm}$.
At $\lambda = \lambda_0$, $\phi_{\alpha,\lambda_0,\pm}$ are colinear (that is, the Jacobian of $\phi_{\alpha,\lambda_0,\pm}$ vanishes). 
Up to a renormalisation we may assume that
$\phi_{\alpha,\lambda_0,+} = \phi_{\alpha,\lambda_0,-}$.  
For $\lambda \ne \lambda_0$ the solution of $Ray_{\alpha,\lambda}(\phi) = \psi$ is explicitely given by
\begin{equation}\label{def-psRphi}
\phi(z) = \phi_{\alpha,\lambda,+}(z) \int_z^{+ \infty} {\phi_{\alpha,\lambda,-}(x) \over Jac(x)} \psi(x) dx
+ \phi_{\alpha,\lambda,-}(z) \int_0^z {\phi_{\alpha,\lambda,+}(x) \over Jac(x)} \psi(x) dx
\end{equation}
where
$$
Jac(x) := \phi_{\alpha,\lambda,-}(x) \partial_x \phi_{\alpha,\lambda,+}(x) 
- \phi_{\alpha,\lambda,+}(x) \partial_x \phi_{\alpha,\lambda,-}(x) 
$$
is the Jacobian of $\phi_{\alpha,\lambda,\pm}$. 
Note that, as $\lambda_0$ is a simple eigenvalue, $Jac(\lambda_0) = 0$ and that for $\lambda$ near $\lambda_0$,
$$
Jac(\lambda) = (\lambda - \lambda_0) \widetilde  Jac(\lambda)
$$
where $\widetilde  Jac(\lambda)$ is a smooth function with $\widetilde  Jac(\lambda_0)\not =0$ since $\lambda_0$
is a simple eigenvalue. Let us also define
$$
\widetilde  \phi_{\alpha,\lambda,\pm} = {\phi_{\alpha,\lambda,\pm} - \phi_{\alpha,\lambda_0,\pm} \over \lambda - \lambda_0} .
$$
Then it follows from \eqref{def-psRphi} that 
$$
\phi(z) = {\phi_{\alpha,\lambda_0,+}(z) \over \lambda - \lambda_0} 
\int_0^{+ \infty} {\phi_{\alpha,\lambda_0,+}(x) \over \widetilde  Jac(x)} \psi(x) dx + \widetilde  \phi(z) 
$$
where 
\begin{equation}\label{def-psphi}
\widetilde  \phi(z) = \int_0^{+ \infty} \widetilde  G(x,z) \psi(x) dx 
\end{equation}
with $$
\widetilde  G(x,z) = {\widetilde  \phi_{\alpha,\lambda,+}(z) \phi_{\alpha,\lambda,-}(x) + \phi_{\alpha,\lambda,+} (z) \widetilde  \phi_{\alpha,\lambda,-}(x) 
+ (\lambda - \lambda_0) \widetilde  \phi_{\alpha,\lambda,+}(z) \widetilde  \phi_{\alpha,\lambda,-}(x) \over \widetilde  Jac(x)}
$$
if $x > z$, and a similar expression if $x < z$.
This computation may be rewritten as follows.
 Let $l$ be the linear form defined by
$$
l(\psi) =  \int_0^{+ \infty} {\phi_{\alpha,\lambda_0,+}(x) \over \widetilde  Jac(x)} \psi(x) dx.
$$
Then, for any $\psi$,  if $l(\psi) = 0$ then $\widetilde  \phi$ solves
$Ray_{\alpha,\lambda}(\widetilde  \phi) = \psi$. 
In particular, as the image of the Rayleigh operator $Im(Ray_{\alpha,\lambda_0})$ is of codimension $1$, 
$Ker(l) = Im(Ray_{\alpha,\lambda_0})$. 
Note that, as $\lambda_0$ is a simple eigenvalue,
$\phi_{\alpha,\lambda_0,+}$ is not in $Im(Ray_{\alpha,\lambda_0})$. Therefore, $l(\phi_{\alpha,\lambda_0,+} ) \ne 0$.
As a consequence
$$
\tilde \psi = \psi - {l(\psi) \over l(\phi_{\alpha,\lambda_0,+} )}  \phi_{\alpha,\lambda_0,+} \in Im(Ray_{\alpha,\lambda})
$$ 
since the image by $l$ of this function vanishes.
We then have
\beq \label{decomp1}
Ray_{\alpha,\lambda}(\widetilde  \phi) = \psi -  {l(\psi) \over l(\phi_{\alpha,\lambda_0,+} )}   \phi_{\alpha,\lambda_0,+},
\eeq
where
$$
\tilde \phi(z) = \int_0^{+ \infty} \tilde G(x,z) \tilde \psi(x) dx .
$$
That is, $\widetilde \phi$ defines the pseudoinverse $Ray_{\alpha,\lambda}^{-1}$ of $Ray_{\alpha,\lambda}$ for $\lambda$ near $\lambda_0$. We shall now fulfill a similar analysis for the $Orr_{\alpha,\lambda}$ operator.


\subsection{Orr Sommerfeld equation}


Let us now prove Theorem \ref{theopseudo}. 
We follow the analysis in the previous section to construct the Green function $\widetilde  G_{\alpha,\lambda}(x,z)$ 
for the pseudoinverse of $Orr_{\alpha,\lambda}$. Let $\lambda_0^{app}$ be a simple eigenvalue of the approximate
Evans function $E^{app}$. $Ray_{\alpha,\lambda}$ operator. To simplify the notation we drop the "app" and set
$\lambda_0 = \lambda_0^{app}$.
At $\lambda = \lambda_0$, the matrix $M$, defined by (\ref{matriceM}), is singular since its first two columns are colinear.
Up to the multiplication by a constant of $\phi_{s,-}$, we may assume that $\phi_{s,\pm}$ coincide at $\lambda = \lambda_0$.
To desingularize it we introduce
$$
\Lambda = \left( \begin{array}{cccc} 
(\lambda - \lambda_0)^{-1} & 1 & 0 & 0 \cr
- (\lambda - \lambda_0)^{-1} & 1 & 0 & 0 \cr
0 & 0 & 1 & 0 \cr
0 & 0 & 0 & 1 \cr
\end{array} \right) .
$$
Then, recalling \eqref{def-vvvapp} and defining $\widetilde  M = M \Lambda$, with the notation of (\ref{Mv}), we have
\begin{equation}\label{def-psv}
v = \Lambda \widetilde  M^{-1} (0,0,0,1/ \nu \mu_f^3) ,
\end{equation}
when $\lambda \ne \lambda_0$.
The arguments applied to the matrix $M$ in Section $3$ may now be applied to $\widetilde  M$ since the corresponding matrix
$$
\widetilde  A_2 = A_2 \Lambda
= \left( \begin{array}{cc} (\phi_{s,-} - \phi_{s,+}) / (\lambda - \lambda_0) & \phi_{s,+} \cr
(\partial_y \phi_{s,-} - \partial_y \phi_{s,+}) / (\lambda - \lambda_0) & \partial_y \phi_{s,+} \cr
\end{array} \right)
$$  
is non singular near $\lambda = \lambda_0$, keeping in mind that $\lambda_0$ is a simple eigenvalue. 

Let $l_4 = (l_{4,1},...,l_{4,4})$ 
be the fourth line of the inverse of $\widetilde  M$. 
It  follows from \eqref{def-psv} that 
$$
v = \Lambda l_4(x) / \nu \mu_f^3.
$$
The singular part $v^s$ of $v$, namely the terms involving $(\lambda - \lambda_0)^{-1}$, is 
$$
v^s ={1 \over \nu \mu_f^3} l_{4,1}(x) (1 , -1, 0, 0).
$$ 
Let us now compute $l_{4,1}(x)$. We have to evaluate $\Lambda^{-1} A_2^{-1} A_1^{-1} B D^{-1}(0,1)$ (see Section \ref{sec-bda}).
But, up to higher order terms, $A_1^{-1} B D^{-1} \sim (0,\mu_f)$.  Note that
$$
A_2^{-1} = {1 \over E^{app}(\alpha,\lambda)} \left( \begin{array}{cc}
\partial_y \phi_{s,+} & - \phi_{s,+} \cr
- \partial_y \phi_{s,-} & \phi_{s,-} \cr \end{array} \right).
$$
Hence when $\lambda$ is close to $\lambda_0$, 
$$
A_2^{-1} A_1^{-1} B D^{-1}(0,1) \sim {\mu_f \over E^{app}(\alpha,\lambda)} \phi_{s,+} (-1,1)
$$ 
namely like $C (-\mu_f,  \mu_f) \phi_{s,+} / (\lambda - \lambda_0)$.
At leading order, the computation is exactly the same as in the previous section.
Let 
$$
L(\psi) = -  \int_0^{+ \infty} l_{4,1}(x) \psi(x) dx .
$$
Then, at leading order, $L = l$. Moreover, the regular part $v^r$ of $v = v^r + v^s$ is
$$
v^r =  {1 \over \nu \mu_f^3}  \Bigl( l_{4,2},l_{4,2},l_{4,3},l_{4,4} \Bigr) .
$$
We now define $\widetilde  G^{app}(x,z)$ to be the approximate Green kernel that corresponds to the regular part $v^r$, 
recalling the Green function construction in \eqref{def-GappX}-\eqref{def-vvvapp}.
Setting 
$$
\widetilde  \psi =  \psi - {L(\psi) \over L(\phi_{\alpha,\lambda_0,+})}  \phi_{\alpha,\lambda_0,+}  ,
$$
we have $L(\widetilde  \psi) = 0$ and so 
$$
Orr_{\alpha,\lambda}(\widetilde  G^{app} \star \psi) = \widetilde  \psi .
$$
The exact Green function $\widetilde G_{\alpha,\lambda}(x,z)$ then follows by iteration
as in Section \ref{sec-exactGreen}. 


\bibliographystyle{abbrv}

\def\cprime{$'$} \def\cprime{$'$}

\end{document}